\documentclass[12pt]{amsart}

\usepackage{amsmath}
\usepackage{latexsym}
\usepackage{amsfonts}
\usepackage{amssymb}

\newtheorem{Thm}{Theorem}
\newtheorem{thm}{Theorem}[section]
\newtheorem{lem}[Thm]{Lemma}

\theoremstyle{definition}

\newcommand{\R}{\mathbb R}

\newcommand{\N}{\mathbb N}
\newcommand{\D}{\mathbb D}

\begin{document}

\title{Poisson type generators for $L^1(\R)$}%

\author{Gerard Ascensi}
\address{Departament de Matem{\`a}tiques, Universitat Aut{\`o}noma de Barcelona, 08193 Cerdanyola del Valles, SPAIN}
\email{gascensi@mat.uab.cat}

\author{Joaquim Bruna}
\address{Departament de Matem{\`a}tiques, Universitat Aut{\`o}noma de Barcelona, 08193 Cerdanyola del Valles, SPAIN}
\email{bruna@mat.uab.cat}

\thanks{Both authors are suported by MTM2005-08984-C02-01 and 2005SGR00611 projects.}

\begin{abstract}
We characterize the discrete sets $\Lambda\subseteq\R$ such that
$\{\varphi(t-\lambda), \lambda\in\Lambda \}$ span $L^1(\R)$,
$\varphi$ being an $L^1(\R)$-function whose Fourier transform
behaves like $e^{-2\pi |\xi|}$.
\end{abstract}

\maketitle

\section{Introduction}

The study of the generators by translations for $L^p(\R)$ has been
a classical topic of study in harmonic analysis. Results in
\cite{Br} and \cite{BOU} characterize the discrete sets
$\Lambda\subseteq\R$ for which there exists a function $\varphi\in
L^1(\R)$ with the property that $\{\varphi(t-\lambda),
\lambda\in\Lambda \}$ span $L^1(\R)$ as those having infinite
Beurling-Malliavin density. In \cite{Ole} and \cite{OlU} we can
find results proving that in $L^2(\R)$ there are more sets whit
this property, and that a characterization in terms of densities
is not posible.

Given a function $\varphi$, a natural problem is to characterize
the discrete sets $\Lambda$ such that $\{\varphi(t-\lambda),
\lambda\in\Lambda \}$ span $L^p(\R)$. There are very few complete
results of this kind. In \cite{BrM} Bruna and Melnikov give a
complete characterization for the Poisson function:

\begin{equation*}
P(t)=\frac{1}{\pi}\frac{1}{1+t^2}.
\end{equation*}

\begin{thm}[Bruna, Melnikov]\label{brunamel}
The translates $\{P(t-\lambda), \lambda\in\Lambda\}$ span
$L^p(\R)$, $1\leqslant p<\infty$ if and only if

\begin{equation}
\sum_{\lambda\in\Lambda}e^{-\frac{\pi}{2}|\lambda|}=\infty.
\end{equation}
\end{thm}

In the $L^2(\R)$ case, by Fourier transform,
$\{\varphi(t-\lambda), \lambda\in\Lambda \}$ span $L^2(\R)$ if and
only if the set of exponentials $\{e^{2\pi i\lambda\xi},
\lambda\in\Lambda\}$ span the weighted space $L^2(\R,
|\widehat{\varphi}|^2)$, and hence the above characterization
holds for any function $\varphi$ such that
$|\widehat{\varphi}|\simeq e^{-2\pi|\xi|}$. The aim of this note
is to give a generalization of this type for the $L^1(\R)$ case.

\begin{thm}\label{teoquasiPoisson}
Assume $\varphi\in L^1(\R)$ has non-vanishing Fourier transform
satisfying
$$A e^{-2\pi|\xi|}\leq |\widehat{\varphi}(\xi)|\leq B
e^{-2\pi|\xi|}$$ for some constants $A,B$. Assume also that
$|(\widehat{\varphi})'|(\xi)=O(e^{-2\pi|\xi|})$. Then
$\{\varphi(t-\lambda), \lambda\in\Lambda \}$ span $L^1(\R)$ if and
only if condition \ref{brunamel} holds.
\end{thm}
In fact the proof will show that condition \ref{brunamel} is
necessary if $Ae^{-2\pi|\xi|}\leqslant|\widehat{\varphi}(\xi)|$
and it is sufficient if both $\widehat{\varphi}(\xi)$ and
$(\widehat{\varphi})'(\xi)$ are $O(e^{-2\pi|\xi|})$.

\section{Proof of the theorem}

Notice first that if certain translates of $\varphi\in L^1(\R)$
span $L^1(\R)$, then obviously $\widehat{\varphi}(\xi)\neq 0$ for
all $\xi$. In fact, non-vanishing of $\widehat{\varphi}$
characterizes (as a consequence of Wiener's Tauberian theorem)
those $\varphi$ such that all its translates span $L^1(\R)$.
Analogously, for $p=2$, a necessary condition in order than some
translates of $\varphi$ span $L^2(\R)$ is that
$\widehat{\varphi}(\xi)\neq 0$ for almost all $\xi$, this being
equivalent to the fact that all translates of $\varphi $ span
$L^2(\R)$.

\begin{lem} Assume $h\in L^1(\R)$ and that $\widehat{h}(\xi)\neq
0$ for all $\xi$. Then, if $f\in L^p(\R), 1\leq p\leq \infty$, and
the convolution $f*h$ is zero, then $f=0$. The same holds if $h\in
L^2(\R)$ and $\widehat{h}(\xi)\neq 0$ almost everywhere.
\end{lem}
\begin{proof}
For $1\leq p\leq 2$, the Fourier transform of $f$ is a function in
$L^{q}, \frac 1p+\frac 1q=1$, and the Fourier transform of $f*h$
is $\widehat{f}\,\widehat{h}$, so the lemma follows. In the
general case, we consider the closed subspace $E$ of $L^1(\R)$
consisting of functions $g$ such that $f*g=0$; since $E$ is
translation invariant and contains $h$, Wiener's tauberian theorem
implies that $E$ is the whole $L^1(\R)$, and this implies $f=0$.
When $h\in L^2(\R)$ we use Beurling's theorem describing all
closed translation-invariant subspaces of $L^2(\R)$ to reach the
same result.
\end{proof}

We note in passing that the version of Wiener or Beurling theorem
for $L^p(\R), 1<p<2$, is still unknown.

In the following, we assume that $\varphi\in L^1(\R)\cap L^p(\R)$
satisfies $\widehat{\varphi}(\xi)\neq 0$ for all $\xi$ or else
$\varphi\in L^2(\R)\cap L^p(\R)$ with $\widehat{\varphi}(\xi)\neq
0$ for almost all $\xi$. By duality, $\{\varphi(t-\lambda),
\lambda\in\Lambda \}$ span $L^p(\R)$ if and only if $f\in L^q(\R)$
and
\begin{equation*}
\widetilde{\varphi}*f(\lambda)=\int_{\R}f(t) \varphi(t-\lambda)
\,dt=0 \quad \forall \lambda\in\Lambda
\end{equation*}
implies $f=0$. Here $\widetilde{\varphi}(t)=\varphi(-t)$. By the
lemma, $f=0$ is equivalent to $\widetilde{\varphi}*f=0$, whence
$\{\varphi(t-\lambda), \lambda\in\Lambda \}$ span $L^p(\R)$ if and
only if $\Lambda$ is a uniqueness set for the space
$$E_{\varphi}^q=\{F=f*\widetilde{\varphi}; f\in L^q(\R)\}$$
meaning that $F\in E_{\varphi}^q, F(\lambda)=0,
\lambda\in\Lambda$, implies $F=0$.

\begin{lem}\label{generaciopertotselsp} Assume $h\in L^1(\R)$ and $\widehat{h}(\xi)\neq 0 $ for
every $\xi$ (respectively, $h\in L^2(\R)$ with
$\widehat{h}(\xi)\neq 0$ almost everywhere). Then if
$\{\varphi(t-\lambda), \lambda\in\Lambda \}$ span $L^1(\R)$ then
$\{(\varphi\ast h)(t-\lambda), \lambda\in\Lambda\}$ span $L^1(\R)$
(respectively $L^2(\R)$).
\end{lem}
\begin{proof} For $f\in L^q(\R)$,
\begin{equation*}
\int_{\R} f(t)\, (\varphi\ast h)(t-\lambda)\,
dt=\int_{\R}(\widetilde{h}\ast f)(x)\varphi(x-\lambda)\, dx,
\end{equation*}
whence the result follows from the lemma above.
\end{proof}

\begin{lem}\label{lemmaA}
If $\phi(t)=P\ast\widehat{P}(t)$  then $\{\phi(t-\lambda),
 \lambda\in\Lambda\}$ span $L^p(\R), 1\leq p<\infty$ if and only if
\ref{brunamel} holds
\end{lem}
\begin{proof}
It is clear that this condition is sufficient, since $\phi$ is a
convolution of $P$ with a function of $L^1(\R)$ and we can apply
lemma \ref{generaciopertotselsp}. For the necessity we will revise
the proof of theorem \ref{brunamel}. By duality, we must see that
if $\sum_{\lambda\in\Lambda}e^{-\frac{\pi}{2}|\lambda|}<\infty$
then we can find $g\in L^q(\R), \neq 0$ such that:
\begin{equation*}
\int_{\R}g(t)\phi(t-\lambda)\,dt=0 \qquad \forall
\lambda\in\Lambda,
\end{equation*}
where we can think that $g$ is real. The above integral equals
\begin{equation*}
\int_{\R}g(t)\bigl(P\ast\widehat{P}\bigr)(t-\lambda)\,dt
=\frac{1}{\pi}\int_{\R}\bigl(g\ast\widehat{P}\bigr)(t)\frac{1}{1+(t-\lambda)^2}\,dt
\end{equation*}
Now we complexify this expression:
\begin{equation}\label{descripFenEB}
F(z)=\frac{1}{\pi}\int_{\R}\frac{f(t)}{(t-z)^2+1}\,dt
\end{equation}
with $f=g\ast\widehat{P}$. When $f$ ranges in $L^q(\R)$, $F$
ranges in the space  $E^q(B)$ which in \cite{BrM} is shown to be
the space of holomorphic functions in $B={\Im z<1}$ such that:
\begin{align*}
&\|F\|= \sup_{|y|<1}\int_{\R}|\Re F(x+iy)|^q\,dx=\|F\|_q^q<\infty\\
&F(\overline{z})=\overline{F(z)},\quad z\in B.
\end{align*}
For $E^{\infty}(B)$  the first condition is replaced by $\Re F$
bounded. Hence we have to find $F\in E^q(B)$ such that
$F(\lambda)=0$ for every $\lambda\in\Lambda$ and that it can be
written as \eqref{descripFenEB} with $f=g\ast\widehat{P}$ for some
$g\in L^q(\R)$. We use the Fourier transform to see that:
\begin{equation*}
F(x)=\int_{\R}\widehat{f}(\xi)e^{-2\pi|\xi|}e^{2\pi i x\xi}\,d\xi.
\end{equation*}
By analytical continuation we obtain:
\begin{equation*}
F(z)=\int_{\R}\widehat{f}(\xi)e^{-2\pi|\xi|}e^{2\pi i z\xi}\,d\xi,
\end{equation*}
where we want that
$\widehat{f}(\xi)=\frac{1}{\pi}\frac{\widehat{g}(\xi)}{1+\xi^2}$
with $g\in L^q(\R)$. That is, we search $F\in E^q(B)$ that can be
written as:
\begin{equation*}
F(z)=\frac{1}{\pi}\int_{\R}\frac{\widehat{g}(\xi)}{1+\xi^2}e^{-2\pi|\xi|}e^{2\pi
i z\xi}\,d\xi
\end{equation*}
with $g\in L^q(\R)$. Since
\begin{equation*}
F''(z)=\int_{\R}\widehat{f}(\xi)(2\pi i\xi)^2e^{-2\pi|\xi|}e^{2\pi i
z\xi}\,d\xi,
\end{equation*}
this amounts to $F''\in E^q(B)$. So we have reduced the problem to
find $F\in E^q(B)$ such that $F(\lambda)=0$ for
$\lambda\in\Lambda$ and such that $F''\in E^q(B)$.

Now, as in \cite{BrM}, we translate the problem to the disk. The
comformal map from $B$ to the disk is given by
\begin{equation*}
w=\Phi(z)=\frac{e^{\frac{\pi}{2}z}-1}{e^{\frac{\pi}{2}z}+1}.
\end{equation*}
Let $\Gamma=\Phi(\Lambda)\subset \R\cap\D$, as shown in \cite{BrM}
finiteness of
$\sum_{\lambda\in\Lambda}e^{-\frac{\pi}{2}|\lambda|}$ is
equivalent to
\begin{equation*}
\sum_{\gamma\in\Gamma}\log \frac{1}{|\gamma|}< \infty,
\end{equation*}
that is the Blaschke condition. This guarantees that the product
\begin{equation*}
\beta(w)=\prod_{\gamma\in\Gamma}\frac{w-\gamma}{1-\gamma w}
\end{equation*}
is convergent (it is necessary to multiply by $w$ if $0\in\Gamma$).

We suppose that $H$ is a holomorphic function of the disk. If
$F(z)=H(\Phi(z))$ then $F\in E^q(B)$ exactly when $g(s)=F(s\pm
i)=H\left(\frac{ie^{\frac{\pi}{2}s}-1}{ie^{\frac{\pi}{2}s}+1}\right)$
is in $L^q(\R)$, that is if
\begin{equation}\label{acotacioE(b)1}
\int_{\R}|g(s)|^q\,ds=\frac{1}{\pi}\int_{|z|=1}|H(z)|^q\frac{|dz|}{|1-z^2|}<\infty.
\end{equation}
We need $F''\in E^q(B)$ as well. Computing
\begin{equation*}
F''(z)=H''\bigl(\Phi(z)\bigr)\left(\frac{\pi
e^{\frac{\pi}{2}z}}{(e^{\frac{\pi}{2}z}+1)^2}\right)^2+
H'\bigl(\Phi(z)\bigr)\frac{\pi^2}{2}\frac{e^{\frac{\pi}{2}z}(1-e^{\frac{\pi}{2}z})}{(e^{\frac{\pi}{2}z}+1)^3}
\end{equation*}
and changing variables again we need
\begin{align}
\int_{|z|=1}\left|H''(z)\frac{\pi i (z+1)(z-1)}{\bigl(i(z+1)+(z-1)\bigr)^2}\right|^q\frac{|dz|}
{|1-z^2|} <\infty\label{acotacioE(b)2}\\
\int_{|z|=1}\left|
H'(z)\frac{\pi^2}{2}\frac{i(z+1)(z-1)\bigl((z-1)-i(z+1)\bigr)}{\bigl(i(z+1)+(z+1)\bigr)^3}\right|^q
\frac{|dz|}{|1-z^2|}
<\infty \label{acotacioE(b)3}.
\end{align}
Therefore we have to find a holomorphic function $H$ in the disk
with $H(\gamma)=0$ for $\gamma\in\Gamma$,
$H(\overline{z})=\overline{H(z)}$ and so that
\eqref{acotacioE(b)1}, \eqref{acotacioE(b)2} and
\eqref{acotacioE(b)3} are fulfilled. We choose
$H(z)=(1-z^2)^n\beta(z)$ with $n$ big enough. We  first bound the
derivatives of $\beta$:
\begin{equation*}
\beta'(z)=\sum_{\gamma\in\Gamma}\frac{1-\gamma^2}{(1-\gamma
z)^2}\prod_{\lambda\in\Gamma, \lambda\neq\gamma}
\frac{z-\lambda}{1-\lambda z}.
\end{equation*}
The product is bounded by $1$ independently of $\gamma$ for almost
all $z$. Moreover, for $|z|=1$ we have  $|1-\gamma
z|=|z-\lambda|\geqslant\frac{1}{2}|1-z^2|$ and
$|1-\gamma^2|\leqslant 2(1-|\gamma|)$. Therefore
\begin{equation*}
|\beta'(z)|\leqslant
\frac{2}{|1-z^2|^2}\sum_{\gamma\in\Gamma}|1-\gamma^2|\leqslant
\frac{4}{|1-z^2|^2}\sum_{\gamma\in\Gamma}1-|\gamma|
\leqslant\frac{2K}{|1-z^2|^2},
\end{equation*}
where we have used the Blaschke condition. For the second
derivative we have
\begin{multline*}
\beta''(z)=2\sum_{\gamma_1\neq\gamma_2\in\Gamma}\frac{1-\gamma_1^2}{(1-\gamma_1
z)^2}
\frac{1-\gamma_2^2}{(1-\gamma_2 z)^2}\prod_{\lambda\in\Gamma, \lambda\neq\gamma_1,\gamma_2}
\frac{z-\lambda}{1-\lambda z}+\\
+2\sum_{\gamma\in\Gamma}\frac{1-\gamma^2}{(1-\gamma
z)^3}\prod_{\lambda\in\Gamma, \lambda\neq\gamma}
\frac{z-\lambda}{1-\lambda z},
\end{multline*}
which similarly can be bound by
\begin{equation*}
|\beta''(z)|\leqslant\frac{12K^2}{|1-z^2|^4}.
\end{equation*}
Then, choosing $H(z)=(1-z^2)^4\beta(z)$ all required conditions
are fulfilled and  the proof is finished.
\end{proof}

\begin{lem}\label{lemmaB}
Let $\psi(t)=P(t)-P''(t)$. Then $\{\psi(t-\lambda),
\lambda\in\Lambda\}$ span $L^p(\R)$ if and only if condition
\ref{brunamel} holds
\end{lem}
\begin{proof}
Notice that $\widehat{\psi}(\xi)=(1+4\pi^2\xi^2)e^{-2\pi|\xi|}$.
The space $E^q_{\psi}$ consists of the functions

\begin{equation*}
F(z)=\langle f(t),\psi(t-z)\rangle=
\int_{\R}\widehat{f}(\xi)(1+4\pi^2\xi^2)e^{-2\pi|\xi|}e^{-2\pi
i\xi z}=G(z)-G''(z)
\end{equation*}
with $G\in E^q(B)$, that is with $f\in L^q(\R)$. Clearly this
space contains $E^q(B)$, and so the condition (1) is necessary,
for if the series converges we already know that there is $H\in
E^q(B)$ vanishing on $\Lambda$.

For the sufficiency we find a growth condition fulfilled by the
second derivative of a function of $E^q(B)$. There is a constant
$c_q$ such that whenever $G$ is holomorphic in a disk $D(a,R)$ of
center $a$ and radius $R$ one has
$$|G''(a)|^q\leq c_q \frac{1}{R^{2q+2}}\int_{D(a,R)} |G(z)|^q
dA(z).$$ Let $G\in E^q(B)$. For $z\in B$ we apply the above to the
ball of center $z$ and radius $\frac{1-|y|}{2}$ ($z=x+iy$) to get
\begin{multline*}
\int_{B}(1-|y|)^{2q}|G''(z)|^q\,dm(z)\\
\leqslant c_q \int_B\frac{1}{(1-|y|)^2}
\int_{B\bigl(z,\frac{1-|y|}{2}\bigr)}|G(w)|^q\,dm(w)\,dm(z).
\end{multline*}
Apply Fubini and noticing that $\bigl\{z:w\in
B\bigl(z,\frac{1-|y|}{2}\bigr)\bigr\} \subset B(w,1-|\Im w|)$ and
that for $z$ in this set $(1-|\Im w|)^2\leq c (1-|y|)^2$ we
obtain
\begin{equation*}
\int_{B}(1-|y|)^{2q}|G''(z)|^q\,dm(z)\leqslant
c_q\int_B|G(w)|^q\,dm(w).
\end{equation*}
This last integral is bounded since
\begin{equation*}
\int_B|G(z)|^q\,dm(z)=\int_{-1}^1\int_{\R}|G(x+iy)|^q\,dx\,dy\leqslant\int_{-1}^1\|G\|^q\,dy=2\|G\|^q.
\end{equation*}
This says that $G''$ is in a Bergman type space. Obviously $G$
satisfies this condition too, and so the above holds with $G$
replaced by $F\in E^q_{\psi}$. We next translate this integral to
the disk.

We can check  that if
$w=\Phi(z)=\frac{e^{\frac{\pi}{2}z}-1}{e^{\frac{\pi}{2}z}+1}$ then
$1-|y|\geqslant 1-|w|$. If  $H(w)=F(\Phi^{-1}(w))$ then
\begin{align*}
\int_{\D}|H(w)|^q(1-|w|)^{2q-1}\,dm(w)\leqslant&\int_{D}(1-|w|)^{2q}|H(w)|^q\frac{dm(w)}{|1-w^2|}\\
\leqslant&\int_{B}(1-|y|)^{2q}|F''(z)|^q\,dm(z).
\end{align*}
Therefore $H$ is in the Bergman space in the disk with weight
$(1-|w|)^{2q-1}$. The set of zeros contained in a diameter of a
function of this space satisfies the Blaschke condition
\cite{Kor}. Therefore the zeros of a function $F\in E^q_{\psi}$
satisfy
$\sum_{\lambda\in\Lambda}e^{-\frac{\pi}{2}|\lambda|}<\infty$ and
so (1) is sufficient.
\end{proof}

Now, theorem \ref{teoquasiPoisson} can be deduced using the
previous lemmas. Assume that $\{\varphi(t-\lambda),
\lambda\in\Lambda \}$ span $L^1(\R)$ and that $
Ae^{-2\pi|\xi|}\leqslant|\widehat{\varphi}(\xi)|$. Writing
$$\phi=P*\widehat{P}=h*\varphi$$
with
$$\widehat{h}(\xi)=\frac{e^{-2\pi|\xi|}}{\widehat{\varphi}(\xi)(1+\xi^2)}.$$
By lemma \ref{generaciopertotselsp} (in the $L^2$ case), the
functions $\{\phi(t-\lambda)\}_{\lambda\in\Lambda}$ span $L^2(\R)$
and therefore by \ref{lemmaA} (1) must hold. Assume next that
$\widehat{\varphi}(\xi)$ and $(\widehat{\varphi})'(\xi)$ are
$O(e^{-2\pi|\xi|})$ and that (1) holds; we write
$$\varphi=h*(P-P'')$$
with
$$\widehat{h}(\xi)=\frac{\widehat{\varphi}(\xi)}{e^{-2\pi|\xi|}(1+4\pi^2\xi^2)}.$$
The hypothesis on $\varphi$ implies that both $\widehat{h}$ and
$(\widehat{h})'$ are in $L^2(\R)$, whence both $h(x)$ and $xh(x)$
are in $L^2(\R)$, so $h\in L^1(\R)$. Since (1) holds, by lemma
\ref{lemmaB} the functions $\{\psi(t-\lambda), \lambda\in\Lambda
\}$ span $L^1(\R)$ and then lemma \ref{generaciopertotselsp}
implies that the functions $\{\varphi(t-\lambda),
\lambda\in\Lambda \}$ span $L^1(\R)$.

\section{Generalizations and comments.}

Using the same ideas we can also prove

\begin{thm}
Fix $n\in\N$. Let $\varphi$ be a function for which there exists
constants $A,B>0$ such that:
\begin{equation}
A\frac{e^{-2\pi|\xi|}}{1+\xi^{2n}}\leqslant|\widehat{\varphi}(\xi)|\leqslant
B(1+\xi^{2n})e^{-2\pi|\xi|}.
\end{equation}
We also suppose that:
\begin{equation*}
|\widehat{\varphi}'(\xi)|\leqslant C(1+\xi^{2n}e^{-2\pi|\xi|}).
\end{equation*}
Then the set $\{\varphi(t-\lambda), \lambda\in\Lambda \}$ spans
$L^1(\R)$ if and only if condition \ref{brunamel} holds.
\end{thm}

Only slight modifications are needed. For instance, one must use
that the $2n$-th derivative of a function of $E^q(B)$ is in a
Bergman type space (whit a different weight); the Korenblum's
result quoted before applies to all the Bergman spaces as it is in
fact true for functions in the class $A^{-\infty}$. Another fact
which is needed is that the $2n$-th derivative of the Blaschke
product appearing above can be bounded by
$\frac{K}{|1-z^2|^{4n}}$.

As used in the proof, for the function  $P$ the space $E^q_P$ is
exactly $E^q(B)$. For the functions $\varphi$ considered here we
have not exactly described this space, yet we can describe its
uniqueness sets.

In \cite{Zal} the Gaussian function $G$ is considered. A complete
characterization is not achieved. In fact, one can show that in
this case the space $E^2_{\varphi}$ can be identified with the
Fock space, for which the description of the uniqueness sets is an
open question. It is known that a sufficient condition in order
that the translates $\{G(t-\lambda), \lambda\in\Lambda\}$ span
$L^2(\R)$ is that the series  $\sum_n
\frac{1}{|\lambda_n|^{2+\varepsilon}}$ diverges for some
$\varepsilon$, while it is necessary that it diverges for
$\varepsilon=0$.

\end{document}